\documentclass{svjour3}                     
\usepackage{graphicx}
\usepackage{pgf,tikz,pgfplots}
\usepackage{tkz-euclide}
\usetikzlibrary{quotes,angles}
\usepackage{amssymb,amsmath}

\begin{document}
	\newtheorem{algo}{Algorithm}
	\newtheorem{proc}{Procedure}
\newtheorem{rem}{Remark}


\title{Deciding whether two quadratic surfaces actually intersect}


\titlerunning{Deciding whether two quadrics intersect}        

\author{Huu-Quang Nguyen       \and
        Ruey-Lin Sheu \and Yong Xia 
}


\institute{Huu-Quang Nguyen \at
              Institute of Natural Science Education, Vinh
              University, Vinh, Nghe An, Vietnam\\Department of Mathematics, National Cheng Kung University, Tainan, Taiwan\\
                           \email{quangdhv@gmail.com}
           \and
           Ruey-Lin Sheu \at
              Department of Mathematics, National Cheng Kung University, Tainan, Taiwan\\
              \email{rsheu@mail.ncku.edu.tw}
              \and
              Yong Xia \at
              LMIB of the Ministry of Education, School of Mathematical Sciences, Beihang University, Beijing 100191, China\\
              \email{dearyxia@gmail.com}
}

\date{Received: date / Accepted: date}

\maketitle

\begin{abstract}
We answer an open question proposed by P\'{o}lik and Terlaky in 2007 that: {\it how we can decide whether two quadratic surfaces intersect without actually computing the intersections?}

\keywords{QCQP, Joint numerical range, S-procedure, Separation of quadrics, Intersection of quadrics, Attainability.}
\end{abstract}
%

\section{Introduction}

Quadratic surfaces (also known as {\it quadrics}) are common in mechanical \textsc{CAD/CAM}. The quadrics intersection problem (QSIC) arises in computer graphics; see Levin \cite[1979]{L79} for a standard reference. It studies the morphology of the intersection curve, and finds ways to compute the mutual curve.
Due to the complication of all possible types of intersections there might arise,
deciding as a priori whether the quadrics actually intersect (the decision version of (QSIC)) helps to render
the intersection in a faster way (Wilf and Manor \cite[1993]{Wilf-Manor}).

The decision version of (QSIC) turns out to be
difficult in its own right. Even for two ellipsoids in the three dimensional space, the classical approaches to determine whether they intersect - either to do it by a constrained
numerical optimization or to rely on numerical observations of the eigenvalues, both are not satisfactory solutions. We need an analytical mathematical solution to answer the intersection of all types, not just for ellipsoids in 3D, but for all cases in general $n$ dimensions (Chan \cite[2006]{Chan06}). It was also posed formally as an open question by P\'{o}lik and Terlaky \cite[2007]{Polik-Terlaky07} that: {\it how we can decide whether two quadratic surfaces intersect without actually computing the intersections?} We shall answer the question in this paper.

Let the two quadratic surfaces be represented by equations $f_i(x)=0, i=1,2,$  where $f_1(x)=x^TA_1x+2a^Tx+a_0, f_2(x)=x^TA_2x+2b^Tx+b_0$ are real-valued quadratic functions of $n$-variables with symmetric matrices $A_1, A_2;$  $a,b\in\mathbb{R}^n$ and $a_0,b_0\in \mathbb{R}.$ Let
$$\mathcal{P}=\{x\in \mathbb{R}^n: f_1(x)=f_2(x)=0\}.$$
The quadrics intersection problem is now the feasibility problem of the set $\mathcal{P}$. To deal with it, we consider a type of polynomial optimization problem of degree 4:
\begin{eqnarray*}
	{\rm (Po4)}\hspace*{2cm}
	\begin{array}{lll}
		v{\rm (Po4)}= &\inf\limits_{x\in \mathbb{R}^n} & (f_1(x))^2+(f_2(x))^2.
	\end{array}
\end{eqnarray*}
It is obvious that $\mathcal{P}= \emptyset$ if and only if either $v(\rm Po4)>0$ or $v(\rm Po4)=0$ but unachieved. The following example shows that the unattainability of (Po4) is possible.
\begin{example}\label{}
Let $f_1(x,y)=x^2, f_2(x,y)=xy-1.$ Consider the sequence
$\{(\frac{1}{n}, n)^T\}_{n=1}^{\infty}$ in $\mathbb{R}^2.$
Since
$$\lim\limits_{n\to \infty}(f_1(x_n,y_n))^2+(f_2(x_n,y_n))^2=\lim\limits_{n\to \infty}\dfrac{1}{n^2}=0,$$
we know $v{\rm (Po4)}=0.$ However, $\mathcal{P}=\{(x,y): x=0, xy-1=0\}=\emptyset$ implies that $v(\rm Po4)=0$ is not attained.
\end{example}

Our approach relies on several known results in literature. 
\begin{itemize}
\item Attainability of quadratic programming subject to one quadratic inequality constraint:
    \begin{eqnarray*}
    {\rm (QP1QC)}\hspace*{0.3cm}
		\begin{array}{lll}
			&\inf\limits_{x\in \mathbb{R}^n} &f_1(x)\\
			&~~{\rm s.t.}&f_2(x)\le 0.
		\end{array}
	\end{eqnarray*}
The problem (QP1QC) can be solved by the classical S-lemma \cite{MJJ,Polik-Terlaky07,Y03}, while the attainability of (QP1QC) was completely characterized with a set of necessary and sufficient conditions by Hsia, Lin and Sheu \cite[Theorem 4.4, Remark 4.5, 2014]{Xia2014}
\item Quadratic programming subject to one quadratic equality constraint:
\begin{eqnarray*}
		{\rm (QP1EQC)}\hspace*{0.3cm}
		\begin{array}{lll}
			&\inf\limits_{x\in \mathbb{R}^n} &f_1(x)\\
			&~~{\rm s.t.}&f_2(x)= 0.
		\end{array}
\end{eqnarray*}
The problem has been completely solved by Xia, Wang and Sheu \cite[2016]{XWS} by their version of S-lemma with equality.
\item Geometrical representation of ``the S-lemma with equality'' and ``the convexity of the joint numerical range of $f_1$ and $f_2.$'' While the former was first characterized analytically by Xia et al. \cite[2016]{XWS} and the latter by Flores-Baz\'{a}n and Opazo \cite[2016]{FO16}, we have provided alternative geometrical version of them by the notation of ``separation of quadratic level sets'' in \cite[2019]{Q-S2019} and \cite[2020]{Q-C-S2020}, respectively. More separation properties about the quadratic level sets shown in \cite[2020]{Q-S2020} (by Nguyen and Sheu) will be also quoted and used.
\end{itemize}

Combining all the above results, our scheme to check the feasibility of the set $\mathcal{P}$ contains the following procedures.

\begin{itemize}
\item[$\bullet$] If $A_1$ and $A_2$ are linearly independent, we develop a new type of $\mathcal{S}$-procedure with $2$ quadratic constraints (Theorem \ref{th1} of this paper). It enables us to compute $v(\rm Po4)$ via an SDP.
\item[$\bullet$] If $A_1=t^*A_2$ (or $A_2=t^*A_1$) with $t^*\in\mathbb{R},$ by elementary analysis, we show that computing $v(\rm Po4)$ can be reduced to solve an {\rm (QP1EQC)}.
\item[$\bullet$] If $v(\rm Po4)=0,$ we show the attainability of (Po4) can be either reduced to checking whether $\{x\in\mathbb{R}^n: f_i(x)=0\}$ separates $\{x\in\mathbb{R}^n: f_j(x)=0\}$ for some $i,j\in\{1,2\}$; or be reduced to check the attainability of some (QP1QC)'s.
\end{itemize}

\noindent{\bf Notations}

Throughout this paper, for a symmetric matrix $A$, $A\succ 0$ $(A\succeq 0)$ denotes $A$
is positive definite (positive semidefinite); $A^{\dagger}$ denotes the Moore-Penrose
generalized inverse of $A$; $\mathcal{R}(A)$
and $\mathcal{N}(A)$ denote its Range and Null spaces, respectively. The notion $\{f \star 0\}$ represents the set $\{x\in \mathbb{R}^n: f(x)\star 0\}$, where $\star\in \{<, \leq, =\}$. The set $\mathbf{R}(f_1,f_2) = \left\{\left(f_1(x),f_2(x)\right)^T|~x \in \mathbb{R}^n \right\}$ is called the joint numerical range of $f_1,f_2.$

\section{Preliminary results}

This section is devoted to a brief sketch of what is known in recognizing whether, given two quadrics, one could separate the other.

Let $f(x)=x^TAx+2a^Tx+a_0, g(x)=x^TBx+2b^Tx+b_0$ be two quadratic functions defined on $\mathbb{R}^n$.
If $f(x)$ is non-constant, after appropriate change of variables (\cite[2019]{Q-S2019}), we may assume that $f(x)$ adopts one of the following five canonical forms:
\begin{eqnarray}
&& -x_1^2-\cdots - x_k^2+\delta (x_{k+1}^2+\cdots +
x_m^2)+\theta; \label{form:1}\\
&&-x_1^2-\cdots - x_k^2+\delta(x_{k+1}^2+\cdots + x_m^2)-\theta';\label{form:2}\\
&&-x_1^2-\cdots - x_k^2+\delta (x_{k+1}^2+\cdots+ x_m^2)+x_{m+1};\label{form:3}\\
&&\hskip8pt x_1^2+\cdots + x_m^2+\delta x_{m+1}+c';\label{form:4}\\
&& \hskip8pt \delta x_1+c', \label{form:5}
\end{eqnarray}
where $\delta\in\{0,1\}, \theta\ge0,~\theta'>0.$
\vskip0.1cm
\noindent{\it Remark:
Without loss of generality, we often simplify $\theta\ge0$ as $\theta\in\{0,1\},$ and represent $\theta'>0$ just by $\theta'=1.$}

\begin{definition}{\rm\cite{Q-S2019,Q-S2020}} \label{def-sep}
	The level set $\{g=0\}$ is said to separate the set $\{f\star 0\}$, where $\star \in \{<,=,\leq\}$, if there are non-empty subsets $L^-$ and $L^+$ of $\{f \star 0\}$ such that
	\begin{equation} \label{eq:def_sep}
	\begin{aligned}
	&L^-\cup L^+ = \{f \star 0\} ~ \text{ and } \\
	&g(u^-)g(u^+)<0,~\forall~ u^-\in L^-;~\forall~u^+\in L^+.
	\end{aligned}
	\end{equation}
\end{definition}

From the definition, we see that, if $\{g=0\}$ separates $\{f\star0\}$, then $g$ satisfies the two-side Slater's condition, i.e., there are $y,z$ such that $g(y)<0<g(z)$. Moreover, from \eqref{eq:def_sep}, $ 0\not\in g(\{f\star0\}).$ Therefore,
      \begin{equation}\label{S-lemma-type}
        \{g=0\} \subset \{f\star0\}^c \mbox{ or equivalently, } \{f\star0\}\cap \{g=0\}=\emptyset.
      \end{equation}
The case that matters us here is to replace ``$\star$'' with ``$=$'' in \eqref{S-lemma-type}. Then, either $\{g=0\}$ separates $\{f=0\}$ or $\{f=0\}$ separates $\{g=0\}$ provides a sufficient condition for null-intersection of them.

\begin{theorem}{\rm \cite[Theorem 1]{Q-S2019}}\label{th3rd} 
    The hypersurface $\{g = 0\}$ separates $\{f < 0\}$ if and only if
    \begin{itemize}
        \item[$\rm(i)$] $f(x)$ is of the form $-x_1^2+\delta(x_2^2+\cdots
        +x_m^2)+\theta,~ \delta\in\{0, 1\},~\theta\ge0;$
        \item[$\rm(ii)$] With the same basis and $\delta$ as in (i),~ $g(x)$ has the affine form $b_1x_1+\delta(b_2x_2+\cdots
        +b_mx_m)+b_0, b_1\ne 0;$
        \item[$\rm(iii)$]
        $f|_{\{g = 0\}}(x)=-(\delta\dfrac{b_2}{b_1}x_2+\cdots+\delta\dfrac{b_m}{b_1}x_m+\dfrac{b_0}{b_1})^2+\delta(x_2^2+\cdots
        +x_m^2)+\theta\geq 0,~\forall \\ (x_2, \cdots,
        x_n)^T\in\mathbb{R}^{n-1}.$ That is, $\{g=0\}\subset\{f\ge0\}.$
    \end{itemize}
\end{theorem}

In \cite[2019]{Q-S2019},  Nguyen and Sheu first gave the above definition and showed that the hypersurface $\{g = 0\}$ separates $\{f < 0\}$ is equivalent to the S-lemma with equality. The S-lemma with equality is concerned with the conditions under which the following two statements are
equivalent:
\begin{itemize}
    \item[](${\rm E_1}$)~~ ($\forall x\in\mathbb R^n$)
    $~g(x)= 0~\Longrightarrow~ f(x)\ge 0.$ 
    \item[](${\rm E_2}$)~~   ($\exists \lambda\in\mathbb R$)
    ($\forall x\in\mathbb R^n$)~$f(x) + \lambda g(x)\ge0.$ 
\end{itemize}
\begin{theorem}{\rm \cite[Theorem 2]{Q-S2019}}\label{th4th} If $g(x)$ takes both positive and negative values,
    then ``$(\rm E_1)$ holds, but $(\rm E_2)$ fails''
    if and only if $\{g=0\}$ separates $\{f <0\}.$
\end{theorem}

Notice that Yakubovich \cite[1971]{Yakubovich71} asserted that, if $g(x)\le0$ satisfies Slater's condition, then (${\rm S_1}$) and (${\rm S_2}$) below are always equivalent without any other conditions.
\begin{itemize}
 \item[](${\rm S_1}$)~~ ($\forall x\in\Bbb R^n$) $~g(x)\le 0~\Longrightarrow~ f(x)\ge 0.$ 
 \item[](${\rm S_2}$)~~   There exists a $\mu\ge0$ such that
$f(x) + \mu g(x)\ge0, ~\forall x\in \mathbb{R}^n.$ 
\end{itemize}
This is the classical S-lemma.

Separation property was also proved to link equivalently to the convexity of the joint numerical range.

\begin{theorem}{\rm \cite[Theorem 3.1]{Q-C-S2020}} \label{thm:noncvx_sep}
The joint numerical range $\mathbf{R}(f,g)$ is non-convex if and only if there exists $\alpha, \beta \in \mathbb{R}$ such that $\{g=\beta\}$ separates $\{f=\alpha\}$ or $\{f=\alpha\}$ separates $\{g=\beta\}$. More precisely,
\begin{enumerate}
\item[(a)] when both $f$ and $g$ are quadratic functions, $\mathbf{R}(f,g)$ is non-convex if and only if $(\exists \alpha, \beta \in \mathbb{R})~\{f=\alpha\}$ and $\{g=\beta\}$ mutually separates each other;
\item[(b)] when one of $f$ and $g$ is affine, $\mathbf{R}(f,g)$ is non-convex if and only if $(\exists \alpha, \beta \in \mathbb{R})$ the affine $\beta$-level set separates the quadratic $\alpha$-level set.
\end{enumerate}
\end{theorem}

The following theorem shows that checking whether $\{g=0\}$ separates $\{f=0\}$ can be reduced to checking whether an affine space separates $\{f=0\}$.
\begin{theorem}{\rm \cite[Theorem 2.4]{Q-S2020}}\label{thm:separation_reduce_linear}
	The set $\{g=0\}$ separates the set $\{f=0\}$ if and only if $B=\lambda A$ for some $\lambda\in\mathbb{R}$ and $\{-\lambda f+g=0\}$ separates $\{f=0\}$.
\end{theorem}

\begin{theorem}{\rm \cite[Lemma 2.5]{Q-S2020}}\label{th8xx}
	Suppose that $g(x)$ is an affine function, then $\{g=0\}$ separates $\{f=0\}$ if and only if there exists a basis such that
	\begin{itemize}
		\item[$\rm(i)$] $f(x)$ has the form $-x_1^2+\delta(x_2^2+\cdots
		+x_m^2)+\theta,~ \delta \in\{0, 1\},~\theta>0;$
		
		\item[$\rm(ii)$] With the same basis, $g(x)$ has the form $b_1x_1+\delta(b_2x_2+\cdots+b_mx_m)+b_0, b_1\ne 0;$
		
		\item[$\rm(iii)$]
		$f(x)|_{\{g = 0\}}=-(\delta\sum_{i=2}^{m}\dfrac{b_i}{b_1}x_i+\dfrac{b_0}{b_1})^2+\delta\sum_{i=2}^{m}x_i^2+1> 0,~\forall (x_2, \cdots,
		x_n)^T\in\mathbb{R}^{n-1}.$
	\end{itemize}
	or
	
	\begin{itemize}
		\item[$\rm(i)'$] $-f(x)$  has the form $-x_1^2+\delta(x_2^2+\cdots
		+x_m^2)+\theta,~ \delta \in\{0, 1\},~\theta>0;$
		
		\item[$\rm(ii)'$] With the same basis,  $g(x)$ has the form $b_1x_1+\delta(b_2x_2+\cdots+b_mx_m)+b_0, b_1\ne 0;$
		
		\item[$\rm(iii)'$]
		$-f(x)|_{\{g = 0\}}=-(\delta\sum_{i=2}^{m}\dfrac{b_i}{b_1}x_i+\dfrac{b_0}{b_1})^2+\delta\sum_{i=2}^{m}x_i^2+1> 0,~\forall (x_2, \cdots,
		x_n)^T\in\mathbb{R}^{n-1}.$
	\end{itemize}
\end{theorem}

\begin{theorem}{\rm \cite[Theorem 4.16]{FO16},\cite[Theorem 4.2]{Q-C-S2020}}\label{th001}
If $\{A, B\}$ are linearly independent, the joint numerical range $\mathbf{R}(f,g)$ is a convex set in $\mathbb{R}^2$.
\end{theorem}

There are some additional properties about separation to be used later.

\begin{proposition}{\rm \cite[Proposition 1]{Q-C-S2020}}\label{prop:linear_comb_preserve_separation}
If $\{g=0\}$ separates $\{f=0\}$, then $\{\eta f+\theta g=0\}$ separates $\{\sigma f=0\}$ for all $\eta, \theta, \sigma \in \mathbb{R}$ with $\theta \neq 0$, $\sigma \neq 0$.
\end{proposition}



%

\section{Computing $v(\rm Po4)$}
\subsection{A new S-procedure}

Let $f_1(x)=x^TA_1x+2a^Tx+a_0,~f_2(x)=x^TA_2x+2b^Tx+b_0$  and $F: \mathbb{R}^2\to \mathbb{R}$ be defined as $F(z)=z_1^2+z_2^2,~z=(z_1,z_2)^T\in\mathbb{R}^2.$ We have the following new type of $\mathcal{S}$-procedure.

\begin{theorem}\label{th1}
	Suppose that $\{A_1, A_2\}$ are linearly independent and $\gamma\in\mathbb{R}.$  The following two statements are equivalent:
	\begin{itemize}
		\item[]${\rm (G_1)}$ $(\forall x\in\mathbb{R}^n,z\in\mathbb{R}^2)$ $f_1(x)-z_1=0, f_2(x)-z_2=0$  $~\Rightarrow~ F(z)-\gamma \geq 0.$ 
		\item[]${\rm (G_2)}$~$(\exists \alpha,\beta\in\mathbb{R})$
		$F(z)-\gamma +\alpha(f_1(x)-z_1)+\beta(f_2(x)-z_2)\geq 0,~~ \forall (x, z)\in\mathbb{R}^n\times\mathbb{R}^{2}$. 
	\end{itemize}
\end{theorem}

\proof{} Note that $\rm(G_2) \Rightarrow \rm(G_1)$ is trivial. We only prove that $\rm(G_1) \Rightarrow \rm(G_2)$.

Since $\{A_1, A_2\}$ are linearly independent, by Theorem \ref{th001},
$$\mathbf{R}(f_1,f_2)=\{z\in \mathbb{R}^2|~z_1=f_1(x),~z_2=f_2(x),~ x\in \mathbb{R}^n\}$$
is convex. Due to $\rm(G_1),$ we have $F(z)-\gamma\geq 0,~ \forall z\in \mathbf{R}(f_1,f_2),$
which implies that
$$\mathbf{R}(f_1,f_2)\cap D=\emptyset,$$ where $D:=\{z\in \mathbb{R}^2|~ F(z)-\gamma<0\}.$
Since $F(z)$ is convex, $D$ is convex and open. By the separating hyperplane theorem, there exists $v=(\bar{\alpha}, \bar{\beta})^T\not=0$ such that
%
\begin{eqnarray}
\bar{\alpha}z_1+\bar{\beta}z_2+\bar{\gamma}\geq 0,&~\forall ~ z\in \mathbf{R}(f_1,f_2), \label{pt4}\\
\bar{\alpha}z_1+\bar{\beta}z_2+\bar{\gamma}< 0,&~ \forall ~z\in D. \label{pt5}
\end{eqnarray}
From (\ref{pt5}), $\bar{\alpha}z_1+\bar{\beta}z_2+\bar{\gamma}\geq 0 \Rightarrow z\not\in D \Rightarrow F(z)-\gamma\geq 0.$  By the classical S-lemma, there exists $t\geq 0$ such that
\begin{equation}\label{pt5a}F(z)-\gamma-t(\bar{\alpha}z_1+\bar{\beta}z_2+\bar{\gamma})\ge 0,~ \forall ~z\in \mathbb{R}^2.\end{equation}

If $t=0$, choose $\alpha=\beta=0$. Then, $\rm(G_2)$ holds.

If $t>0$, by (\ref{pt4}), one has, $\forall~x\in \mathbb{R}^n,$
\begin{eqnarray}
&&t\bar{\alpha}f_1(x)+t\bar{\beta}f_2(x)+t\bar{\gamma}\geq 0 \nonumber\\ 
&\Leftrightarrow& t\bar{\alpha} f_1(x) +t\bar{\beta} f_2(x)+(t\bar{\alpha}-t\bar{\alpha})z_1+(t\bar{\beta}-t\bar{\beta})z_2+ t\bar{\gamma}\geq 0 \nonumber \\ 
&\Leftrightarrow& t\bar{\alpha}(f_1(x)-z_1)+t\bar{\beta}(f_2(x)-z_2) +(t\bar{\alpha}z_1+t\bar{\beta}z_2+t\bar{\gamma})\ge0. \label{pt5f}
\end{eqnarray}
Let $\alpha=t\bar{\alpha},~\beta=t\bar{\beta}.$ Then, (\ref{pt5f}) renders $\rm(G_2)$ by
(\ref{pt5a}). The proof is complete. \hfill$\Box$
\endproof

The following example shows that linear independence of $\{A_1, A_2\}$ is necessary for Theorem \ref{th1}.

\begin{example} Let $f_1(x)=x_1,~f_2(x)=-x_1^2+x_2^2+1,$ $~F(z)=z_1^2+z_2^2$ and $\gamma=0.25.$ Then, $A_1=0$ and $\{A_1,A_2\}$ are thus linearly dependent.
We assert that $\rm(G_1)$ holds for this example. That is,
\begin{equation}\label{ex2a}
F(z)-0.25=x_1^2+(-x_1^2+x_2^2+1)^2-0.25\geq 0,~\forall x_1, x_2\in\mathbb{R}.
\end{equation}
Note that \eqref{ex2a} is obviously true for $|x_1|\ge0.5,~x_2\in\mathbb{R}.$
When $|x_1|<0.5,$
$$x_1^2+(-x_1^2+x_2^2+1)^2-0.25\ge(-x_1^2+x_2^2+1)^2-0.25>(-0.25+1)^2-0.25>0.$$
Therefore, $\rm(G_1)$ holds. On the other hand, to validate $\rm(G_2)$
	\begin{eqnarray*}
		&& F(z)-0.25+\alpha(f_1(x)-z_1)+\beta(f_2(x)-z_2) \\
		&=&z_1^2+z_2^2-0.25+\alpha(x_1-z_1)+\beta(-x_1^2+x_2^2+1-z_2)\geq 0,
	\end{eqnarray*}
it is necessary that, for some choice of $\alpha,\beta,$
	\begin{equation*}\label{inequa1}M=\begin{pmatrix}
	-\beta & ~~0 &~ 0 &0&\alpha/2\\
	0 & ~~\beta & 0 & 0&0\\
	0 & ~~0 & 1 & 0&-\alpha/2\\
	0 & ~~0 & 0 & 1&-\beta/2\\
	\alpha/2 &~~ 0 & -\alpha/2 & -\beta/2&\beta-0.25
	\end{pmatrix}\succeq 0.\end{equation*}
	It implies that $\beta=\alpha=0,$ but then
    $$\begin{pmatrix}
	-\beta&~~\alpha/2\\
	\alpha/2&~\beta-0.25
	\end{pmatrix}=\begin{pmatrix}
	0&~~0\\
	0&~-0.25
	\end{pmatrix}\not\succeq0.$$ So, $\rm(G_2)$ must fail. \end{example}

%

\subsection{Solving the optimal value of (Po4)}\label{sec22}

Let us first write the problem (Po4) as a quadratic programming with two quadratic equality constraints as follows.
	\begin{eqnarray}\label{ap:0100}
	\begin{array}{lll}
	&\inf\limits_{(x, z)\in \mathbb{R}^{n}\times \mathbb{R}^{2}} & (z_1)^2+(z_2)^2\\
	&\hspace*{1cm}{\rm s.t.}& \begin{cases}\begin{array}{ll}
	f_1(x)-z_1=0, \label{nls:001} \\
	f_2(x)-z_2=0. \label{nls:002}
	\end{array}\end{cases}
	\end{array}
	\end{eqnarray}

\noindent {\bf Case 1}: $\{A_1, A_2\}$ are linearly independent.
Applying the new $\mathcal{S}$-procedure in Theorem \ref{th1}, we show that the optimal value $v(\rm{Po4})$  can be obtained by solving an SDP.

\begin{theorem}\label{thm3asdf} Under the condition $\{A_1, A_2\}$ are linearly independent, the optimal value of $\rm (Po4)$, $v(\rm{Po4}),$ can be computed by
	\begin{align}\label{19aa}
	v(\rm{Po4})=\underset{\tiny \begin{array}{lll}
		\gamma, \,\alpha, \,\beta \in \mathbb{R}
		\end{array}}{\sup}  \left\lbrace\begin{array}{cc} \gamma&|~M \succeq 0 \end{array}\right\rbrace,
	\end{align}
	where $M\in\mathbb{R}^{(n+3)\times(n+3)}$ is
	\begin{equation}\label{sigma.1}
	\begin{pmatrix}\begin{array}{ccc}
	1 & \quad \quad ~& 0 \\ 0 &\quad & 1
	\end{array} & ~~[0]_{2\times n} & \begin{array}{cc}
	\frac{-\alpha}{2}  \\ \frac{-\beta}{2}
	\end{array}\\  ~~~[0]_{n\times 2} & \quad \alpha A_1+\beta A_2 & \alpha a+\beta b\\ \begin{array}{cc}
	\frac{-\alpha}{2}  & \quad \frac{-\beta}{2}
	\end{array} & \quad ~\alpha a^T+\beta b^T & \quad \alpha a_0+\beta b_0-\gamma\end{pmatrix},
	\end{equation}where
	 $[0]=(0_{ij})_{2\times n}$ with $0_{ij}=0 ~\forall i, j.$ \end{theorem}

\proof
Since (Po4) can be formulated as \eqref{nls:001}, we have
%
%

\begin{equation}
\begin{array}{lll}
v{\rm (Po4)}&=&\inf\limits_{(x, z)\in \mathbb{R}^n\times \mathbb{R}^2} F(z) \\
& ~~~{\rm s.t.} & \begin{cases}\begin{array}{lll}
f_1(x)-z_1=0 \\
f_2(x)-z_2=0
\end{array}\end{cases}\\
&=&{\sup} \left\lbrace \gamma \left|~ \left\lbrace
(x, z)\in \mathbb{R}^n\times \mathbb{R}^2 \left|~ \begin{array}{ll}
F(z) < \gamma\\
f_1(x)=z_1\\
f_2(x)=z_2
\end{array}\right.
\right\rbrace=\emptyset\right.\right\rbrace\\
&=&\underset{\tiny\begin{array}{lll}
	\gamma, \,\alpha, \,\beta \in \mathbb{R}
	\end{array}}{\sup} \big\{ \gamma |~\overline{F}(x,z)\geq 0~\forall (x, z)\in \mathbb{R}^n\times \mathbb{R}^2\big\},
\end{array}\label{19uio}
\end{equation}
where $\overline{F}(x, z)=F(z)-\gamma +\alpha(f_1(x)-z_1)+\beta(f_2(x)-z_2)$, and the last equality in (\ref{19uio}) holds by Theorem \ref{th1}. Note that $\overline{F}(x,z)\ge0$ can be written as a linear matrix inequality \eqref{19aa} with a matrix $M$ defined in \eqref{sigma.1}. \hfill$\Box$  \endproof
\vskip0.2cm
{\bf Case 2}: Suppose that $\{A_1, A_2\}$ are linearly dependent, say $A_2=t^*A_1.$ If $A_1=A_2=0,$ then (Po4) is an unconstrained convex quadratic optimization problem, which can be solved directly. Hence, assume that $A_1\not=0.$ Then, multiplying the first equation in \eqref{nls:001} by $t^*$ and subtract it from the second equation, we obtain
\begin{eqnarray}
\begin{cases}\begin{array}{lll}
x^TA_1x+a^Tx+a_0-z_1=0,\\
(b^T-t^*a^T)x+(b_0-t^*a_0)+t^*z_1-z_2=0.\label{27uhhu}
\end{array}\end{cases}
\end{eqnarray}
Define $y^T=[x^T,~{z^T}]$; $\bar{a}^T=[a^T, -1, 0],$ $h^T=[b^T-t^*a^T, t^*, -1]$,~$h_0=b_0-t^*a_0$;
and
$$\bar{A}=\begin{pmatrix}\begin{array}{ll}[0]_{n\times n} & [0]_{n\times 2} \\
[0]_{2\times n} &I_{2\times 2}\end{array}\end{pmatrix},~ \bar{A_1}=\begin{pmatrix}\begin{array}{ll}A_1 & [0]_{n\times 2} \\
[0]_{2\times n} &[0]_{2\times 2}\end{array}\end{pmatrix}\not=0.$$
Then, \eqref{ap:0100} becomes
\begin{eqnarray}
&\inf_{y\in \mathbb{R}^{n+2}}& y^T\bar{A}y  \nonumber\\
&{\rm s.t.}&y^T\bar{A_1}y+\bar{a}^Ty+a_0=0,~h^Ty+h_0=0.\nonumber
\end{eqnarray}
By the null space representation, one can write the hyperplane $h^Ty+h_0=0$
as $y=y_0+Vz$ where $V\in\mathbb{R }^{(n+2) \times (n+1)}$ is the matrix basis of  $\mathcal {N}(h),~z\in \mathbb{R}^{n+1}$ and $y_{0} = -\frac{h_0}{h^{T}h}h$.
Then, \eqref{ap:0100} is reduced to the following (QP1EQC) problem:
\begin{eqnarray}
&\inf\limits_{z\in \mathbb{R}^{n+1}}&~ (y_0+Vz)^T\bar{A}(y_0+Vz)\hskip 1.5cm \label{dfgh}\\
&{\rm s.t.}&~(y_0+Vz)^T\bar{A_1}(y_0+Vz)+\bar{a}^T(y_0+Vz)+a_0=0,\nonumber
\end{eqnarray}
which can be solved by applying the S-lemma with equality. See \cite{XWS}.

%
%



\section{Determining the feasibility of the set $\mathcal{P}$}

Suppose the optimal value of problem (Po4) has been computed. If $v(\rm Po4)>0$, we know immediately that the two quadrics $\{f_1=0\}$ and $\{f_2=0\}$ have no intersection. Otherwise, $v(\rm Po4)=0.$ It remains to study whether $v(\rm Po4)=0$ is attainable or not. Our analysis relies on whether the two hypersurfaces $\{f_1=0\}$ and $\{f_2=0\}$ ``separate'' one by the other. See Definition \ref{def-sep} in Section 2. To this end,
we need to assume that both $f_1(x), f_2(x)$ satisfy two-side Slater's condition, i.e., there are $x, y\in \mathbb{R}^n$ such $f_i(x)<0<f_i(y), i\in\{1,2\}$.


If not, let $f_1(x)\geq 0,~\forall x\in\mathbb{R}^n.$ Namely, $f_1(x)$ is a convex function. Then $A_1\succeq 0, a\in\mathcal{R}(A_1), a_0-a^TA_1^+a\geq 0.$ It indicates that the first order equation $A_1x+a=0$ has a solution, say $x_0=-A_1^+a,$ and $f(x_0)=-a^TA^+_1a+a_0=\inf\{f(x):~x\in\mathbb{R}^n\}.$

If $a_0-a^TA^+a> 0,$ then  $\{f_1=0\}=\emptyset.$ In this case, the two quadrics $\{f_1=0\}$ and $\{f_2=0\}$ do not intersect.

Otherwise, $a_0-a^TA^+a=0.$ Then, $\{f_1=0\}=\{x: Ax+a=0\}=\{-A_1^+a+Zy: y\in\mathbb{R}^m\}$ is affine in $\mathbb{R}^n$, where $Z\in\mathbb{ R }^{n\times m}$ is a matrix basis of $\mathcal{N}(A_1)$. Therefore,
\begin{eqnarray*}
&&\{f_1=0\}\cap\{f_2=0\}
=\emptyset \\ &\Leftrightarrow& \{y:(-A_1^+a+Zy)^TA_2(-A_1^+a+Zy)+2b^T(-A_1^+a+Zy)+b_0=0\}=\emptyset.
\end{eqnarray*}
Let $C=Z^TA_2Z, c=-a^TA_1^+A_2Z+b^TZ, c_0=a^TA_1^+A_1^+a-2b^TA_1^+a+b_0.$
Then, the two quadrics do not intersect if,
either $C\succeq 0, c\in\mathcal{R}(C)$ and $c_0-c^TCc>0;$ or $C\preceq 0, c\in\mathcal{R}(C)$ and $c_0-c^TCc<0.$

In order to get the main result of this section, we need the following lemma.

\begin{lemma}\label{l1} Assume that $f_1(x), f_2(x)$ both satisfy two-side Slater's condition; $v(\rm Po4)=0$; $\mathcal{P}=\emptyset$ and there doesn't exist $i\ne j$ in $\{1, 2\}$ such that $\{f_i=0\}$ separates $\{f_j=0\}$. Then, the joint range set $\mathbf{R}(f_1,f_2)=\{(f_1(x), f_2(x))^T: x\in\mathbb{R}^n\}$ is convex.
\end{lemma}
\proof
Suppose on the contrary that $\mathbf{R}(f_1,f_2)$ is not convex. By Theorem \ref{thm:noncvx_sep} in Section 2, there are  $\alpha, \beta\in \mathbb{R}$, $i\ne j$ in $\{1, 2\}$, say $i=1, j=2$ such that $\{f_1=\alpha\}$ separates $\{f_2=\beta\}$.
By Theorem \ref{thm:separation_reduce_linear}, there is some $\lambda\in\mathbb{R}$ such that $A_1=\lambda A_2$ and
{\begin{equation}\label{14poi}
\text{the hyperplane } \{f_1-\lambda f_2-\alpha+\lambda\beta=0\} \text{ separates } \{f_2-\beta=0\}.
\end{equation}
By Theorem \ref{th8xx}, there exists a basis of $\mathbb{R}^n$ such that $f_2-\beta$ or $-(f_2-\beta)$ has the form $$-x_1^2+\delta(x_2^2+\cdots +x_m^2)+c',~ \text{ where } \delta\in\{0,1\},~c'>0$$ and
\begin{equation}\label{mki}
f_1-\lambda f_2-\alpha+\lambda\beta=\nu_1x_1+\delta(\nu_2x_2+\cdots+\nu_mx_m)+\nu_0, ~ \nu_1\ne 0.
\end{equation}
%
%
%
Without loss of generality, we may assume that $$f_2(x)-\beta=-x_1^2+\delta(x_2^2+\cdots +x_m^2)+c',~c'>0$$ since
$-(f_2(x)-\beta)=-x_1^2+\delta(x_2^2+\cdots +x_m^2)+c'$ can be similarly proved.
Due to $\mathcal{P}=\emptyset$, i.e. $\{f_1=0\}\cap\{f_2=0\}=\emptyset$, it implies immediately that}
\begin{equation}\label{f112}
\{f_1-\lambda f_2=0\}\cap\{f_2=0\}=\emptyset.
\end{equation}

\noindent $\bullet$ {Case 1: ${\rm rank}(A_2)\geq 2$}. That is, $\delta=1,~m\ge2.$
We first claim that, from \eqref{f112} there indeed is
\begin{equation}\label{f122}\{f_1-\lambda f_2=0\}\subset \{f_2>0\} \text{   or  } \{f_1-\lambda f_2=0\}\subset \{f_2<0\}.\end{equation}
Indeed, by \eqref{f112}, if \eqref{f122} failed, there would be $\bar{x}, \hat{x}\in \{f_1-\lambda f_2=0\}$ such that $f_2(\bar{x})<0<f_2(\hat{x})$.
From \eqref{mki}, $\{f_1-\lambda f_2=0\}$ is connected, so there exists $\tilde x\in\{f_1-\lambda f_2=0\}$ such that $f_2(\tilde x)=0,$ a contradiction to \eqref{f112}. Therefore, \eqref{f122} must hold true.
	
- If $\{f_1-\lambda f_2=0\}\subset \{f_2>0\},$ namely, $f_2(x)|_{\{f_1-\lambda f_2=0\}}>0,$ the first statement of S-lemma with equality, $(\rm E_1)$, holds. Due to
\begin{equation*}
f_2(x)=-x_1^2+x_2^2+\cdots +x_m^2+c'+\beta,~c>0;
\end{equation*}
and
$f_1-\lambda f_2$ is non-constant affine (see \eqref{mki}), there is no $\gamma\in\mathbb{R}$ such that $f_1(x)-\lambda f_2(x)+\gamma f_2(x)\geq 0, ~\forall x\in\mathbb{R}^n.$ The second statement of S-lemma with equality, $(\rm E_2)$, fails. By Theorem \ref{th4th} and Theorem \ref{th3rd}, we know $\{f_1-\lambda f_2=0\} \text{ separates }  \{f_2<0\},$ and $c'+\beta\ge0.$ In fact, there is $c'+\beta>0.$ Otherwise, as we can see that the following system of equations
\[
\left\{\begin{array}{l}
f_2(x)=-x_1^2+x_2^2+\cdots +x_m^2=0\\
\left(f_1-\lambda f_2\right)(x)=\nu_1x_1+\nu_2x_2+\cdots+\nu_mx_m+\nu_0+\alpha-\lambda\beta=0,~\nu_1\not=0
\end{array}\right.
\]
has a solution
 \[
 \left\{\begin{array}{lr}
x_1=x_2=\frac{-\nu_0-\alpha+\lambda\beta}{\nu_1+\nu_2},x_3=\cdots=x_n=0, & \text{if } \nu_1+\nu_2\not=0~\\
x_1=-x_2=\frac{-\nu_0-\alpha+\lambda\beta}{\nu_1-\nu_2},x_3=\cdots=x_n=0, & \text{if } \nu_1-\nu_2\not=0,
\end{array}\right.
  \]
which is a contradiction to \eqref{f112}.

By $c'+\beta>0;$ $f_1-\lambda f_2$ is non-constant affine; $f_2(x)|_{\{f_1-\lambda f_2=0\}}>0,$ we conclude, by Theorem \ref{th8xx},
that $\{f_1-\lambda f_2=0\}$ separates $\{f_2=0\}.$ It follows from Proposition \ref{prop:linear_comb_preserve_separation} that $\{f_1=0\}$ separates $\{f_2=0\},$
which contradicts to the hypothesis of the lemma. 

%
%

%
%
%


- If $\{f_1-\lambda f_2=0\}\subset \{f_2<0\},$ then  $\{f_1-\lambda f_2=0\}\subset \{-f_2>0\}.$ By the same above arguments for the pair $\{f_1-\lambda f_2, -f_2\}$, we also arrive a contradiction. In a short summary, Case 1 cannot happen.

\noindent $\bullet$ Case 2: ${\rm rank}(A_2) =1$.  Then $\pm (f_2-\beta)$ has the form $-x_1^2+1$, $f_1-\lambda f_2-\alpha+\lambda\beta$ has the form $\nu_1x_1+\nu_0, \nu_1\ne 0$ so that  $\mathbf{R}(f_1-\lambda f_2,f_2)=\{(\nu_1x_1+\nu_0+\alpha-\lambda\beta, \pm(-x_1^2+1)+\beta)^T: x_1\in\mathbb{R}\}$ is a parabola in $\mathbb{R}^2.$ By \eqref{f112}, this parabola does not contain $0$. Let $L_{\lambda}: \mathbb{R}^2 \rightarrow \mathbb{R}^2, L_{\lambda}(u, v)=(u+\lambda v, v)$ be a bijective linear transformation, which maps a parabola to a parabola. Therefore,
$$\mathbf{R}(f_1,f_2)=L_{\lambda}(\mathbf{R}(f_1-\lambda f_2,f_2))$$
is also a parabola that does not contain $0$. It follows that $v({\rm Po4})=\inf\{f_1^2+f_2^2\}>0,$ which is a contradiction. Therefore, Case 2 cannot happen either.

In conclusion, under the hypothesis of the lemma, the joint range set $\mathbf{R}(f_1,f_2)$ must be convex.
 \hfill$\Box$  \endproof

\begin{theorem}\label{thm4} Assume that $f_1(x), f_2(x)$ both satisfy two-side Slater's condition and $v(\rm Po4)=0$. Then $\mathcal{P}=\emptyset$ if and only if either one of the two cases happen.
\begin{itemize}
  \item[(a)] $\{f_i=0\}$ separates $\{f_j=0\}$ for some $i\ne j$ in $\{1, 2\}$;
  \item[(b)] there exists $i\ne j\in\{1, 2\}$ such that at least one of the four (QP1QC)'s happened: $\inf\{\pm f_i(x): \pm f_j(x)\le 0\}=0$. Moreover, for the above (QP1QC)'s that do happen, they must be unattainable.
\end{itemize}
\end{theorem}

\proof $(\Rightarrow)$ If (a) does not hold, we have to prove that (b) is true.
By Lemma \eqref{l1}, the joint range set $\mathbf{R}(f_1,f_2)=\{(f_1(x), f_2(x))^T: x\in\mathbb{R}^n\}$ is convex. Moreover, the range set does not contain the origin due to $\mathcal{P}=\emptyset$. By the separation hyperplane theorem, there exists $(\mu_1, \mu_2)\not=(0,0)$ such that
\begin{equation}\label{13qa}
\mu_1f_1(x)+\mu_2f_2(x)\geq 0, ~\forall x\in\mathbb{R}^n.
\end{equation}
In fact, one has $\mu_1\ne 0$, $\mu_2\ne 0.$ Otherwise, say $\mu_1= 0,\mu_2\ne 0.$ Then $\mu_2f_2(x)\geq 0, ~\forall x\in\mathbb{R}^n$ would have failed the two-side Slater's condition of $f_2(x)$. Now, dividing $|\mu_1|$ at both sides of \eqref{13qa}, we have
\begin{equation}\label{15qa}
sign(\mu_1) f_1(x)+\dfrac{|\mu_2|}{|\mu_1|}(sign(\mu_2)f_2(x))\geq 0, ~\forall x\in\mathbb{R}^n.
\end{equation}
Applying the classical S-lemma to \eqref{15qa}, $\{sign(\mu_2)f_2\leq 0\}\subset \{sign(\mu_1)f_1\geq 0\}$. Or equivalently,
\begin{equation*}\label{16qa}
\{sign(\mu_2)f_2\leq 0\}\cap \{sign(\mu_1)f_1<0\}=\emptyset.
\end{equation*}
Then,
\begin{eqnarray*}\label{17qa}
&&\inf\{sign(\mu_1)f_1(x): x\in \{sign(\mu_2)f_2\leq 0\}\}\\
&=&\sup\{\gamma: \{sign(\mu_1)f_1<\gamma\}\cap\{sign(\mu_2)f_2\leq 0\}=\emptyset\}\\
&\ge&0.
\end{eqnarray*}
We are going to claim that, indeed,
\begin{equation}\label{19qa}
\inf\{sign(\mu_1)f_1(x): x\in \{sign(\mu_2)f_2\leq 0\}\}= 0.
\end{equation}
%
%

Suppose contrarily that $\inf\{sign(\mu_1)f_1(x): x\in \{sign(\mu_2)f_2\leq 0\}\}= \epsilon>0.$ Then,
$\{sign(\mu_1)f_1<\epsilon/2\}\cap\{sign(\mu_2)f_2\leq 0\}=\emptyset.$
By S-lemma, there is $\lambda\geq 0$ such that $sign(\mu_1)f_1(x)-\epsilon/2+\lambda \left(sign(\mu_2)f_2(x)\right)\geq 0,~\forall x$. Hence,
\begin{eqnarray*}
0<\epsilon/2 &\leq&  |sign(\mu_1)f_1(x)+\lambda \left(sign(\mu_2)f_2(x)\right)| \\
 &\leq& \sqrt{1+\lambda^2}\sqrt{(f_1(x))^2+(f_2(x))^2},~\forall x\in\mathbb{R}^n,
\end{eqnarray*}
which leads to $v({\rm Po4})=\inf\{(f_1(x))^2+(f_2(x))^2\}\geq \frac{\epsilon^2}{4(1+\lambda^2)} >0$, a contradiction to our assumption.
So we have proved \eqref{19qa} that one of the four (QP1QC)'s happened: $\inf\{\pm f_i(x): \pm f_j(x)\le 0\}=0$.

Finally, for the (QP1QC) that happened, say $\inf\{f_1(x): f_2(x)\le 0\}=0$, it must be unattainable. Otherwise, let $\hat x$ be an optimal solution to $\inf\{f_1(x): f_2(x)\le 0\}=0$ such that $f_2(\hat x)\le0$ and $f_1(\hat x)=0.$ Since $\mathcal{P}=\emptyset,$ there must be $f_2(\hat x)<0.$ Then, $\inf\{f_1(x): f_2(x)\le 0\}=0$ attains its minimum value $0$ at an interior point. It could happen only when $f_1(x)$ is convex and $f(x)\ge0,~\forall x\in\mathbb{R}^n.$
It violates the assumption that $f_1$ satisfies the two-side Slater's condition.
Therefore, $\inf\{f_1(x): f_2(x)\le 0\}=0$ is unattainable.

$(\Leftarrow)$ It is trivial. Indeed,

\noindent - If (a) holds, by Definition \ref{def-sep}, we immediately have $\{f_1=0\}\cap\{f_2=0\}=\emptyset.$

\noindent - If (b) holds, without loss of generality, we assume that $\inf\{f_1(x): f_2(x)\leq 0\}=0.$ Suppose $x^*\in\{f_1=0\}\cap\{f_2=0\}\not=\emptyset,$ $x^*$ must be a minimizer for $\inf\{f_1(x): f_2(x)\leq 0\}=0,$ contradicting to the assumption that $\inf\{f_1(x): f_2(x)\leq 0\}=0$ is not attainable. Hence, $\{f_1=0\}\cap\{f_2=0\}=\emptyset.$

 \hfill$\Box$  \endproof

Theorem \ref{thm4} reduces the task of checking the quadric intersection problem  to first compute the optimal value of (Po4) as discussed in Subsection \ref{sec22};  to check the separation property of two quadrics as in \cite{Q-S2019,Q-S2020,Q-C-S2020} and to examine the attainability of (QP1QC), $\inf\{\pm f_i(x): \pm f_j(x)\leq  0\}$, as probed in \cite{Xia2014}.
All necessary steps are implementable and they involve only solving an SDP and matrix computation, though the entire procedure is cumbersome.

Examples below are meant to show the validity of our theory.

%
%
%
%
%
%

\begin{example}\label{}
Let $f_1(x)=x_1$, $f_2(x)=x_1x_2-1.$ Then, $\mathcal{P}=\{f_1=0\}\cap\{f_2=0\}=\emptyset.$
By our scheme, we first compute $v({\rm Po4})=\inf\{x_1^2+(x_1x_2-1)^2\}=0.$
It can be observed from the sequence $\{(\frac{1}{n},n)\}\subset \mathbb{R}^2$ on which $\lim\limits_{n \rightarrow\infty}(1/n)^2+(n/n-1)^2=0.$
We now show that $\{f_1=0\}$ separates $\{f_2=0\}$. We also decide that $\mathcal{P}=\{f_1=0\}\cap\{f_2=0\}=\emptyset.$ By change of coordinate with
$(x_1,x_2)=(t+s,t-s),$ $$f_1(t,s)=t+s,~-f_2(t,s)=-t^2+s^2+1 \hbox{ and } -f_2(t,s)|_{t=-s}=1>0.$$
By Theorem \ref{th8xx}, $\{f_1=0\}$ separates $\{f_2=0\}$, which verifies Theorem
\ref{thm4}.
\end{example}

\begin{example}\label{}
Let $f_1(x)=x_1^2+x_1x_2-1$, $f_2(x)=x_1x_2-1$. It is easy to see that $\mathcal{P}=\{f_1=0\}\cap\{f_2=0\}=\emptyset.$ On the other hand, by the sequence $\{(\frac{1}{n},n)\}\subset \mathbb{R}^2,$ since $\lim\limits_{n \rightarrow\infty}(1/n+n/n-1)^2+(n/n-1)^2=0,$ we also have $v({\rm Po4})=0.$

Notice that the two matrices
$$ A_1= \left[\begin{array}{cc}
1 & 0.5\\
0.5 & 0
\end{array}\right];~A_2= \left[\begin{array}{cc}
0 & 0.5\\
0.5 & 0
\end{array}\right]$$
are linearly independent, by Theorem \ref{thm:separation_reduce_linear},
there doesn't occur that $\{f_i=0\}$ separates $\{f_j=0\}$ for any $(i,j)$ pair.
Regarding the 4 (QP1QC)'s $\inf\{\pm f_i(x): \pm f_j(x)\leq  0\}$, there
are only two of them which have the optimal value zero:
\begin{eqnarray}\label{}
\inf\{f_1(x):  -f_2\leq 0\}&=&\inf\{x_1^2+x_1x_2-1: x_1x_2-1\geq 0\}=0; ~~~~\label{21qa}\\
\inf\{-f_2(x):  f_1\leq 0\}&=&\sup\{x_1x_2-1: x_1^2+x_1x_2-1\leq 0\}=0.~~~\label{22qa}
\end{eqnarray}
Both \eqref{21qa} and \eqref{22qa} can be checked to be unattainable. Theorem
\ref{thm4} is thus verified.	\end{example}


\section{Conclusion and Discussion}

In this paper, we propose a scheme to answer conclusively whether two given quadrics intersect or not. We describe the problem as a special type of (unconstrained) polynomial optimization problem of degree 4, which we call it (Po4). The quadric intersection problem is then transformed to solve (Po4) and to detect whether its optimal value $v({\rm Po4})$ is attainable or not. Technically, we develop a new S-procedure to compute $v({\rm Po4})$, while having to
borrow a relatively new concept about the separation of quadratic level sets to conquer the attainability issue of it. We agree that our method still requires further computational details to make it practical, for example, how $v({\rm Po4})=0$ can be rigidly distinguished from $v({\rm Po4})>0$ using today's finite-digit precision machine. However, our purpose is mainly to show that there is a nice, mathematically sensible, approach to the open problem whose effective solution methods were not known to people before. On the other hand,
we feel that the composition of ``quadratic with quadratics'' occurs naturally in many applications, such as those of nonlinear least squares. We hope that our approach can be extended to solving optimization problems involving other quartic polynomials.

\section*{Acknowledgements}

Huu-Quang, Nguyen's research work
was sponsored partially by Taiwan MOST 107-2811-M-006-535 and
Ruey-Lin Sheu's research work
was sponsored partially by Taiwan MOST 107-2115-M-006-011-MY2.

\noindent Xia's research was supported by National Natural Science Foundation of China
under grants 11822103, 11571029, 11771056, and Beijing Natural Science Foundation Z180005.

\bibliographystyle{amsplain}

\begin{thebibliography}{10}

\bibitem {Chan06}
	Chan, K., 2006. \textit{A simple mathematical approach for determining intersection of quadratic surfaces.} In book: Multiscale Optimization Methods and Applications, edited by William W. Hager,~Shu-Jen Huang,~Panos M. Pardalos,~and Oleg A. Prokopyev, pp. 271--298.
	
	
	\bibitem {FO16}
	Flores-Baz\'{a}n, F. and Opazo, F., 2016. \textit{Characterizing the convexity of joint-range for a pair of inhomogeneous quadratic functions and strong duality.} Minimax Theory and its Applications, 1, pp. 257--290.
	
%
%
	
\bibitem {L79}
	Levin, J., 1979. \textit{Mathematical models for determining the intersections of quadratic surfaces.} Computer Graphics and Image Processing, 11, pp. 73--87.
	
\bibitem {MJJ}
	Mor\'{e}, J.J., 1993. \textit{Generalizations of the trust region problem.} Optimization methods and Software, 2(3-4), pp. 189--209.
	
\bibitem{Q-S2019}
Nguyen H.Q. and Sheu, R. L., 2019.  \textit{Geometric properties for level sets of quadratic functions.} Journal of Global Optimization, {73} (2019), 349--369.

\bibitem{Q-S2020}
Nguyen H.Q. and Sheu, R. L., 2020. \textit{Separation -Level Sets.}   Available from: {https://doi.org/10.13140/RG.2.2.25970.12488}.
\bibitem{Q-C-S2020}
Nguyen H.Q., Y.C. Chu and Sheu, R. L. 2020. \textit{On the convexity for the ragne set of two quadratic functions}. Journal of Industrial and Management Optimization (2020, Accepted). Available from: {doi: 10.3934/jimo.2020169}.
	
\bibitem {Polik-Terlaky07}
	Polik, I. and Terlaky, T., 2006. \textit{A survey of the S-lemma.} SIAM Review,  49, pp. 371--418.
	
%
%
	
%
%
%

	\bibitem{Xia2014}
Hsia, Y. and Lin, G. X. and Sheu, R. L., 2014. {\em A revisit to quadratic programming with one inequality quadratic constraint via matrix pencil}. Pacific Journal of Optimization, 10, pp.  461--481.
	
%
%
	
\bibitem {Wilf-Manor}
	Wilf, I. and Manor, Y., 1993. \textit{Quadric-surface intersection curves: shape and structure.} Computer-Aided Design, 25, pp. 633--643.

%
	
\bibitem {XWS}
	Xia, Y. and Wang, S. and Sheu, R.L., 2016. \textit{S-lemma with equality and its applications.} Mathematical Programming Series A, 156(1), pp. 513--547.
	
\bibitem{Yakubovich71}
	Yakubovich, V.A., 1971. {\em S-procedure in nonlinear control theory}. Vestnik Leningrad University (in Russian), 1, pp.  62--77.
	
\bibitem {Y03}
	Ye, Y. and Zhang, S. , 2003. \textit{New results on quadratic minimization.} SIAM Journal on Optimization, 14, pp. 245--267.
	
\end{thebibliography}

\end{document}